elucier@ieee.org

# Fermat's Last Theorem, Solution Sets v7

Ernest R. Lucier

**Abstract**

The non-zero integer solution set is derived for $C^n = A^n + B^n$. The non-zero integer solution set for n = 2 is $[C - (a + b)]^2 = 2ab$. The variables a and b equal $(C - A)$ and $(C - B)$ respectively and are nonzero integer factors of $2M^2$ where M is a non-zero integer. C is greater than $(a + b)$ since the square root of 2ab is an imaginary number when C is less than $(a + b)$. C is equal or greater than $(a + b + 1)$ since we are only considering whole numbers. The derivation of the solution set for n = 2 is applied to n = 3, n = 4, and generalized to n. The solution set for n = n is

$[C - (a + b)]^n = ab([^n_2]C^{n-2}(2) - [^n_3]C^{n-3}(3a + 3b) + \ldots \pm [^n_n] \{[^n_1]a^{n-2} + [^n_2]a^{n-3}b^1 + [^n_3]a^{n-4}b^2 + \ldots + [^n_{n-1}]b^{n-2}\})$. Where the binomial coefficient $[^n_r] = n!/[(n - r)!r!]$ is the coefficient of the $x^r$ term in the polynomial expansion of the binomial power $(1 + x)^n$ and $[^n_r] = 0$ if r > n. Divide this equation by $[C - (a + b)]^{n-2}$ to obtain $[C - (a + b)]^2$. The solution set for $[C - (a + b)]^2$ equals 2ab. The nth solution set equals 2ab only when n equals 2.

$[C - (a + b)]^2$ (I.e., $[C - (a + b)]^n$ divided by $[C - (a + b)]^{n-2}$) is always greater than 2ab when n is greater than 2. Non-zero integer solutions exist only for n = 2.

**Derive a solution for $C^n = A^n + B^n$ with n = 2.**

$C^n = A^n + B^n$ with n = 2 is:

| | |
|---|---|
| $C^2 = A^2 + B^2$ | (1) |

Rewriting A and B in terms of C, a, and b:

| | |
|---|---|
| $A = C - a$ | (2) |
| $B = C - b$ | (3) |

Substitute for A and B in Equation 1:

| | |
|---|---|
| $C^2 = (C - a)^2 + (C - b)^2$ | (4) |

Expand terms

$$C^2 = C^2 - 2Ca + a^2 + C^2 - 2Cb + b^2 \qquad (5)$$

Combine like terms

$$C^2 - 2C(a + b) + a^2 + b^2 = 0 \qquad (6)$$

Complete the square

$$C^2 - 2C(a + b) + (a + b)^2 = \{(a + b)^2 - (a^2 + b^2)\} \qquad (7)$$

and

$$C^2 - 2C(a + b) + (a + b)^2 = 2ab \qquad (8)$$

and

$$[C - (a + b)]^2 = 2ab \qquad (9)$$

Take the square root of both sides of Equation 9:

$$C - (a + b) = \pm (2ab)^{1/2} \qquad (10)$$

Solve for C

$$C = a + b \pm (2ab)^{1/2} \qquad (11)$$

Equation 11 is the solution for $C^n = A^n + B^n$ with $n = 2$.

**Define the non-zero integer solution set for n = 2.**

Factor both sides of Equation 9:

| | |
|---|---|
| $[C - (a + b)]*[C - (a + b)] = 2*a*b$ | (12) |

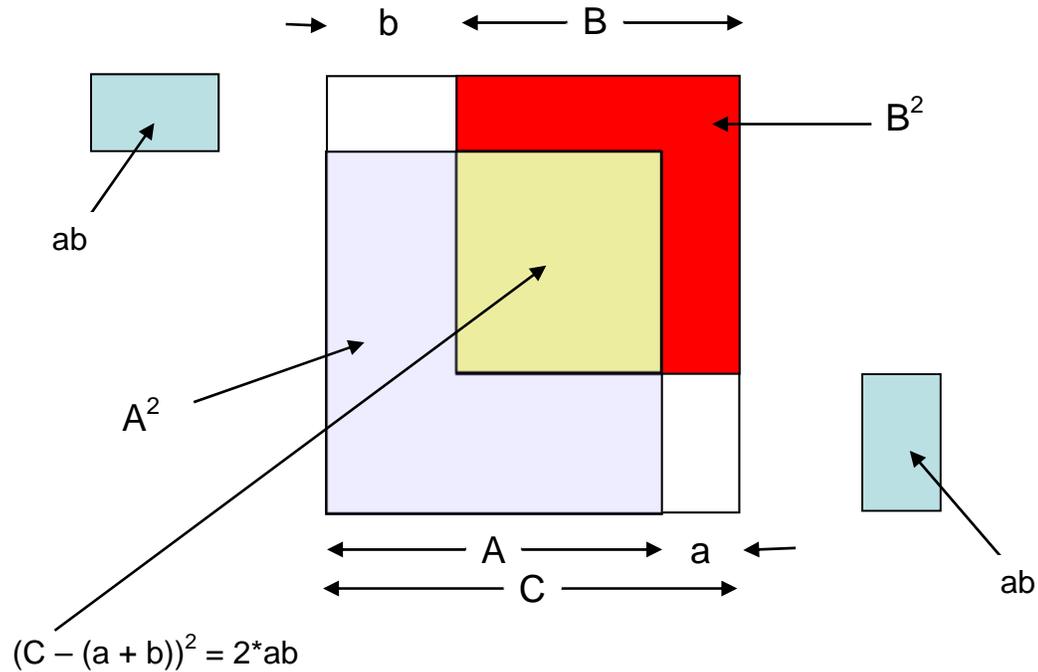

$(C - (a + b))^2 = 2*ab$

Figure 1. Physical representation of Equations 9 and 11 (Not to scale)

Notes:

1. 2*ab completes the square, Equation 9. The area completing the square equals $[C - (a + b)]^2$.

2. $[C - (a + b)]^2$ is area common to both $A^2$ and $B^2$, i.e., overlapping areas. The area inside $C^2$ and outside $A^2$ and $B^2$ (i.e., 2*ab) must equal the common area $[C - (a + b)]^2$ (i.e., common volume). A common area can exist only if $C > (a + b)$.

Define M in terms of a and b to ensure that (2*a*b) of Equation 11 will be a whole number.

Let $2*M^2 = a*b$.

| | |
|---|---|
| $2M^2 = ab$ | (13) |
| $[C - (a + b)]*[C - (a + b)] = 2*2*M^**M$ | (14) |
| Then | |
| $C = a + b \pm 2*M$ | (15) |

Equation 15 is the non-zero integer solution set for Equation 1, $C^2 = A^2 + B^2$, when M is a non-zero integer and a and b are non-zero integer factors of $2M^2$. There are an infinite number of non-zero integer solutions (i.e., Pythagorean triples) for $C^2 = A^2 + B^2$ where M = 1, 2, 3 …

**Example**: let M = 1 so that $2M^2 = 2$

Using Equations 13, 15, 2, and 3:

The factors of $2M^2$ are a = 1 and b = 2 since $2M^2 = 2$.

| | |
|---|---|
| C = 1 + 2 + 2, C = 5, and $C^2 = 25$ | (16) |
| A = C – a, A = 5 – 1, A = 4, and $A^2 = 16$ | (17) |
| B = C – b, B = 5 – 2, B = 3, and $B^2 = 9$ | (18) |

Therefore, when n = 2 and M = 1, this Equation results in a 3, 4, 5 triangle, i.e., 9 + 16 = 25.

**Derive a solution for $C^n = A^n + B^n$ with n = 3.**

$C^n = A^n + B^n$ with n = 3 is:

$$C^3 = A^3 + B^3 \qquad (19)$$

Rewriting A and B in terms of C, a, and b as in n = 2 (Equations 2 and 3):

$$A = C - a \qquad (20)$$
$$B = C - b \qquad (21)$$

Substitute for A and B in Equation 19:

$$C^3 = (C - a)^3 + (C - b)^3 \qquad (22)$$

Expand terms

$$C^3 = C^3 - 3C^2 a + 3Ca^2 - a^3 + C^3 - 3C^2 b + 3Cb^2 - b^3 \qquad (23)$$

Combine like terms

$$C^3 - 3C^2(a + b) + 3C(a^2 + b^2) - (a^3 + b^3) = 0 \qquad (24)$$

Complete the cube

$$C^3 - 3C^2(a + b) + 3C(a + b)^2 - (a + b)^3 = \qquad (25)$$
$$3C\{(a + b)^2 - (a^2 + b^2)\} - \{(a + b)^3 - (a^3 + b^3)\}$$

and

$$C^3 - 3C^2(a + b) + 3C(a + b)^2 - (a + b)^3 = 3C(2ab) - (3a^2 b + 3ab^2) \qquad (26)$$

and

$$[C - (a + b)]^3 = 3ab[2C - (a + b)] \qquad (27)$$

Take the cube root of both sides of Equation 27.

$$C - (a + b) = \{3ab[2C - (a + b)]\}^{1/3} \qquad (28)$$

Solve for C

| | |
|---|---|
| $$C = a + b + \{3ab[2C - (a + b)]\}^{1/3}$$ | (29) |

Equation 29 is the solution for $C^n = A^n + B^n$ with $n = 3$.

**Define the non-zero integer solution set for n = 3.**

Factor both sides of Equation 27:

| | |
|---|---|
| $[C - (a + b)]*[C - (a + b)]*[C - (a + b)] = 3*a*b*[2C - (a + b)]$ | (30) |

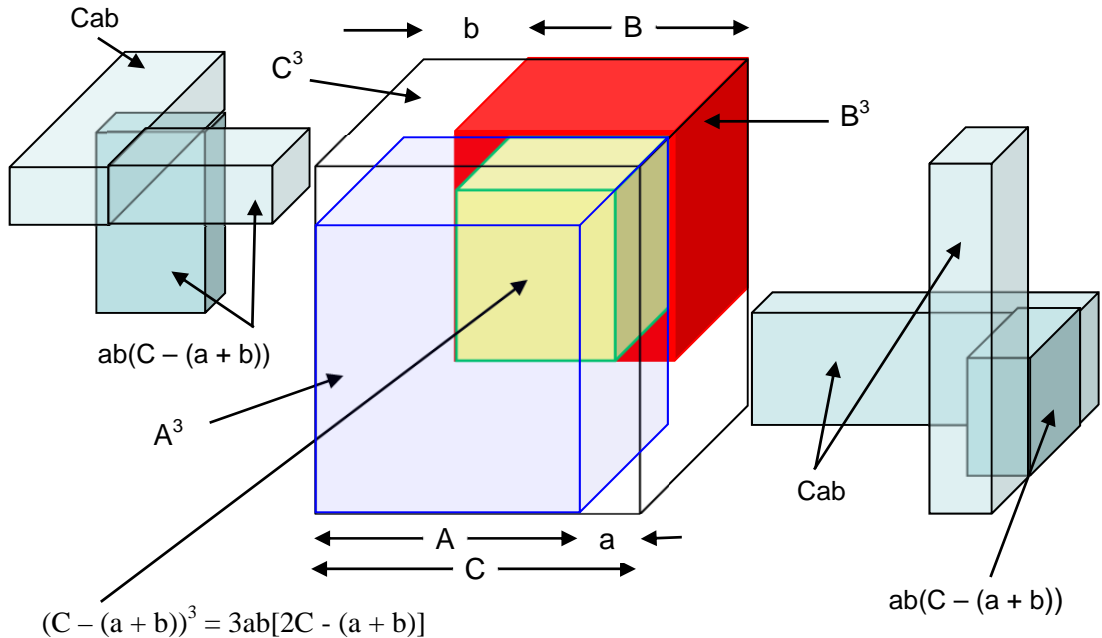

Figure 2. Physical representation of Equation 30 (Not to scale)

Notes:

1. $3*Cab + 3*ab(C - (a + b)) = 3ab[2C - (a + b)]$ completes the cube, Equation 27. The volume completing the cube equals the common volume, $[C - (a + b)]^3$.

2. $[C - (a + b)]^3$ is volume common to both $A^3$ and $B^3$. The volume inside $C^3$ and outside $A^3$ and $B^3$ (i.e., $3ab[2C - (a + b)]$) must equal the common volume $[C - (a + b)]^3$ (i.e., common volume). A common volume exist only if $C > (a + b)$.

Divide both sides of Equation 27 by [C - (a + b)].

| | |
|---|---|
| $[C - (a + b)]*[C - (a + b)] = 3*a*b*[2C - (a + b)] \div [C - (a + b)]$ | (31) |
| $[C - (a + b)]*[C - (a + b)] = 3*a*b*\{C + [C - (a + b)]\} \div [C - (a + b)]$ | (32) |
| $[C - (a + b)]*[C - (a + b)] = 3*a*b*\{C \div [C - (a + b)] + [C - (a + b)] \div [C - (a + b)]\}$ | (33) |
| $[C - (a + b)]*[C - (a + b)] = 3*a*b*\{C \div [C - (a + b)] + 1\}$ | (34) |

Add and subtract (a + b) to C in the numerator on the right side of the Equation 34 and rewrite:

| | |
|---|---|
| $[C - (a + b)]*[C - (a + b)] = 3*a*b*\{[C - (a + b) + (a + b)] \div [C - (a + b)] + 1\}$ | (35) |
| $[C - (a + b)]*[C - (a + b)] = 3*a*b*\{1 + \{(a + b) \div [C - (a + b)]\} + 1\}$ | (36) |
| $[C - (a + b)]*[C - (a + b)] = 3*a*b*\{2 + \{(a + b) \div [C - (a + b)]\}\}$ | (37) |

The non-zero integer solution set for $[C - (a + b)]^2$ from Equation 9 is [2ab] when M is a non-zero integer and a and b are non-zero integer factors of $2M^2$ with $2M^2$ = ab (see Equation 15). Find the value of C in Equation 37 when $[C - (a + b)]^2$ equals [2ab].

Let $[C - (a + b)]*[C - (a + b)]$ equal [2ab] in Equation 37 and solve for C and (a + b).

| | |
|---|---|
| $[2ab] = 3*a*b*\{2 + \{(a + b) \div [C - (a + b)]\}\}$ | (38) |

Divide Equation 38 by 3ab and combine terms:

| | |
|---|---|
| $2/3 = 2 + \{(a + b) \div [C - (a + b)]\}$ | (39) |
| $2/3 - 2 = (a + b) \div [C - (a + b)]$ | (40) |
| $-4/3 = (a + b) \div [C - (a + b)]$ | (41) |
| $- 4/3 * [C - (a + b)] = (a + b)$ | (42) |
| $C - (a + b) = - ¾ * (a + b)$ | (43) |
| $C = - ¾ * (a + b) + (a + b)$ | (44) |
| $C = (1 - ¾) * (a + b)$ | (45) |
| $C = 1/4 * (a + b)$ | (46) |

C is less than (a + b). A common volume exists only when C is greater than (a + b). Therefore, there is no common volume and there are no non-zero integer solutions for n = 3.

**Derive a solution for $C^n = A^n + B^n$ with n = 4.**

$C^n = A^n + B^n$ with n = 4 is:

$$C^4 = A^4 + B^4 \qquad (47)$$

Rewriting A and B in terms of C, a, and b as in n = 2 (Equations 2 and 3):

$$A = C - a \qquad (48)$$
$$B = C - b \qquad (49)$$

Substitute for A and B in Equation 47:

$$C^4 = (C - a)^4 + (C - b)^4 \qquad (50)$$

Expand terms

$$C^4 = C^4 - 4C^3a + 6C^2a^2 - 4Ca^3 + a^4 + C^4 - 4C^3b + 6C^2b^2 - 4Cb^3 + b^4 \qquad (51)$$

Combine like terms

$$C^4 - 4C^3(a + b) + 6C^2(a^2 + b^2) - 4C(a^3 + b^3) + (a^4 + b^4) = 0 \qquad (52)$$

Complete the fourth power

$$C^4 - 4C^3(a+b) + 6C^2(a+b)^2 - 4C(a+b)^3 + (a^4+b^4) =$$
$$6C^2\{(a+b)^2 - (a^2+b^2)\} - 4C\{(a+b)^3 - (a^3+b^3)\} + \{(a+b)^4 - (a^4+b^4)\} \quad (53)$$

and

$$C^4 - 4C^3(a+b) + 6C^2(a+b)^2 - 4C(a+b)^3 + (a^4+b^4) =$$
$$6C^2(2ab) - 4C(3a^2b + 3ab^2) + (4a^3b + 6a^2b^2 + 4ab^3) \quad (54)$$

and

$$[C - (a+b)]^4 = 2ab[6C^2 - 6C(a+b) + (2a^2 + 3ab + 2b^2)] \quad (55)$$

Take the fourth root of both sides of Equation 55.

$$C - (a+b) = \{2ab[6C^2 - 6C(a+b) + (2a^2 + 3ab + 2b^2)]\}^{1/4} \quad (56)$$

Solve for C

$$C = a + b + \{2ab[6C^2 - 6C(a+b) + (2a^2 + 3ab + 2b^2)]\}^{1/4} \quad (57)$$

Equation 57 is the solution for $C^n = A^n + B^n$ with $n = 4$.

**Define the non-zero integer solution set for n = 4.**

Factor both sides of Equation 55.

$$[C - (a + b)]*[C - (a + b)]*[C - (a + b)]*[C - (a + b)]$$
$$= 2ab[6C^2 - 6C(a + b) + (2a^2 + 3ab + 2b^2)] \quad (58)$$

The fourth root of the left side of Equation 58 is $[C - (a + b)]$. The fourth root of the right side must also equal $[C - (a + b)]$ and both sides can be divided by $[C - (a + b)]^2$.

$$[C - (a + b)]^2 = 2ab[6C^2 - 6C(a + b) + (2a^2 + 3ab + 2b^2)] \div [C - (a + b)]^2 \quad (59)$$

C must be greater than $(a + b)$ for a common volume to exist. Let $C = (a + b + x)$ and $x = 1$, i.e., $[C - (a + b)]^2$ equals one (1). Substitute in the denominator of the right side of Equation 59 and solve for $[C - (a + b)]^2$.

$$[C - (a + b)]^2 = 2ab[6C^2 - 6C(a + b) + (2a^2 + 3ab + 2b^2)] \quad (60)$$
$$[C - (a + b)]^2 = 2ab[6(a + b + 1)^2 - 6(a + b + 1)(a + b) + (2a^2 + 3ab + 2b^2)] \quad (61)$$
$$[C - (a + b)]^2 = 2ab*[6(a + b + 1) + (2a^2 + 3ab + 2b^2)] \quad (62)$$

$[C - (a + b)]^2$ must equal $[2ab]$ since the non-zero integer solution set for $[C - (a + b)]*[C - (a + b)]$ from Equation 9 is $[2ab]$ when M is a non-zero integer and a and b are non-zero integer factors of $2M^2$ with $2M^2 = ab$ (see Equation 15).

From Equation 62 $[C - (a + b)]^2$ does not equal $[2ab]$ when C is greater than $(a + b)$ and $n = 4$ since:

$$[2ab] < 2ab*[6(a + b + 1) + (2a^2 + 3ab + 2b^2)] \quad (63)$$
$$\text{for all C greater than } (a + b)$$

Therefore, there is no non-zero integer solution set for n = 4.

**Derive a solution for $C^n = A^n + B^n$ for all n.**

$C^n = A^n + B^n$ for all n:

$$C^n = A^n + B^n \qquad (64)$$

Rewriting A and B in terms of C, a, and b as in n = 2 (Equations 2 and 3):

$$A = C - a \qquad (65)$$
$$B = C - b \qquad (66)$$

Substitute for A and B in Equation 64:

$$C^n = (C - a)^n + (C - b)^n \qquad (67)$$

Expand terms

$$C^n =$$
$$C^n - nC^{n-1}a + [^n_2]C^{n-2}a^2 - [^n_3]C^{n-3}a^3 + \ldots \pm [^n_n]a^n +$$
$$C^n - nC^{n-1}b + [^n_2]C^{n-2}b^2 - [^n_3]C^{n-3}b^3 + \ldots \pm [^n_n]b^n \qquad (68)$$

Where the binomial coefficient $[^n_r] = n!/[(n - r)!r!]$ is the coefficient of the $x^r$ term in the polynomial expansion of the binomial power $(1 + x)^n$ and $[^n_r] = 0$ if r > n.

Combine like terms

$$C^n - nC^{n-1}(a + b) + [^n_2]C^{n-2}(a^2 + b^2) - [^n_3]C^{n-3}(a^3 + b^3) + \ldots \pm [^n_n](a^n + b^n) = 0 \qquad (69)$$

Complete the $n^{th}$ power

$$C^n - nC^{n-1}(a+b) + [^n_2]C^{n-2}(a+b)^2 - [^n_3]C^{n-3}(a+b)^3 + \ldots \pm (a+b)^n =$$
$$[^n_2]C^{n-2}\{(a+b)^2 - (a^2+b^2)\} - [^n_3]C^{n-3}\{(a+b)^3 - (a^3+b^3)\} + \ldots$$
$$\pm [^n_n]C^{n-n}\{(a+b)^n - (a^n+b^n)\} \tag{70}$$

and

$$C^n - nC^{n-1}(a+b) + [^n_2]C^{n-2}(a+b)^2 - [^n_3]C^{n-3}(a+b)^3 + \ldots \pm (a+b)^n =$$
$$[^n_2]C^{n-2}(2ab) - [^n_3]C^{n-3}(3a^2b + 3ab^2) + \ldots$$
$$\pm [^n_n]C^{n-n}\{[^n_1]a^{n-1}b + [^n_2]a^{n-2}b^2 + [^n_3]a^{n-3}b^3 + \ldots + [^n_{n-1}]ab^{n-1}\} \tag{71}$$

and

$$[C - (a+b)]^n =$$
$$[^n_2]C^{n-2}(2ab) - [^n_3]C^{n-3}(3a^2b + 3ab^2) + \ldots$$
$$\pm [^n_n]\{[^n_1]a^{n-1}b + [^n_2]a^{n-2}b^2 + [^n_3]a^{n-3}b^3 + \ldots + [^n_{n-1}]ab^{n-1}\} \tag{72}$$

or

$$[C - (a+b)]^n = ab([^n_2]C^{n-2}(2) - [^n_3]C^{n-3}(3a + 3b) + \ldots$$
$$\pm [^n_n]\{[^n_1]a^{n-2} + [^n_2]a^{n-3}b^1 + [^n_3]a^{n-4}b^2 + \ldots + [^n_{n-1}]b^{n-2}\}) \tag{73}$$

Note: The $[^n_n]$ term is the last term in the series. That is, the first term is the last term when $n = 2$, the second term when $n = 3$, etc.

Take the $n^{th}$ root of both sides of Equation 73.

$$C - (a+b) = a^{1/n}b^{1/n}([^n_2]C^{n-2}(2) - [^n_3]C^{n-3}(3a + 3b) + \ldots$$
$$\pm [^n_n]\{[^n_1]a^{n-2} + [^n_2]a^{n-3}b^1 + [^n_3]a^{n-4}b^2 + \ldots + [^n_{n-1}]b^{n-2}\})^{1/n} \tag{74}$$

Solve for C

$$C = a + b + a^{1/n}b^{1/n}([^n_2]C^{n-2}(2) - [^n_3]C^{n-3}(3a + 3b) + \ldots$$
$$\pm [^n_n]\{[^n_1]a^{n-2} + [^n_2]a^{n-3}b^1 + [^n_3]a^{n-4}b^2 + \ldots + [^n_{n-1}]b^{n-2}\})^{1/n} \tag{75}$$

Equation 75 is the solution for $C^n = A^n + B^n$ for any n.

**Define the non-zero integer solution set for any n.**

Factor both sides of Equation 73.

| | |
|---|---|
| $[C - (a + b)]_1 * [C - (a + b)]_2 * \ldots * [C - (a + b)]_{n-1} * [C - (a + b)]_n =$ <br><br> $a*b*([^n_2]C^{n-2}(2) - [^n_3]C^{n-3}(3a + 3b) + \ldots$ <br><br> $\pm [^n_n] \{[^n_1]a^{n-2} + [^n_2]a^{n-3}b^1 + [^n_3]a^{n-4}b2 + \ldots + [^n_{n-1}]b^{n-2}\})$ | (76) |

The $n^{th}$ root of the left side of Equation 76 is $[C - (a + b)]$. $[C - (a + b)]^n$ can be divided by $[C - (a + b)]^{n-2}$ to obtain $[C - (a + b)]^2$. The non-zero integer solution set for $[C - (a + b)]^2$ is 2ab when a and b are non-zero integer factors of $2M^2$ and M is a non-zero integer. (See Equation 15).

Divide both sides of Equation 73 by $[C - (a + b)]^{n-2}$.

| | |
|---|---|
| $[C - (a + b)]^2 =$ <br><br> $a*b*([^n_2]C^{n-2}(2) - [^n_3]C^{n-3}(3a + 3b) + \ldots$ <br><br> $\pm [^n_n] \{[^n_1]a^{n-2} + [^n_2]a^{n-3}b^1 + [^n_3]a^{n-4}b2 + \ldots + [^n_{n-1}]b^{n-2}\}) \div [C - (a + b)]^{n-2}$ | (77) |

Solve for $[C - (a + b)]^2$ for n = 2, 3, 4, and all n.

For n = 2 rewrite Equation 77 with n = 2.

(Note: The $[^n_n]$ term in Equation 77 is the first term when n = 2):

$$[C - (a + b)]^2 = a*b*([^2_2]C^{2-2}(2)) \div [C - (a + b)]^{2-2} \quad (78)$$

$$[C - (a + b)]^2 = 2ab \quad (79)$$

Equation 79 is the same as Equation 9 (n = 2). The non-zero integer solution set

for $[C - (a + b)]^2$ from Equation 9 is [2ab] when M is a non-zero integer and a and b are non-zero

integer factors of $2M^2$ with $2M^2$ = ab (see Equation 15). There are an infinite number of non-zero

integer solutions (i.e., Pythagorean triples) for $C^2 = A^2 + B^2$ where M = 1, 2, 3 …

For n = 3 rewrite Equation 77 with n = 3.

(Note: The $[^n_n]$ term in Equation 77 is the second term when n = 3):

$$[C - (a + b)]^2 = a*b*([^3_2]C^{3-2}(2) - [^3_3]C^{3-3}(3a + 3b)) \div [C - (a + b)]^{3-2} \quad (80)$$

$$[C - (a + b)]^2 = 3ab(2C - (a + b)) \div [C - (a + b)] \quad (81)$$

C must be greater than (a + b) for a common volume to exist (see Figures 1 and 2). The non-zero integer solution set for $[C - (a + b)]^2$ from Equation 9 is [2ab] when M is a non-zero integer and a and b are non-zero integer factors of $2M^2$ with $2M^2$ = ab (see Equation 15).

From Equation 46, C = 1/4 * (a + b) when n = 3, C is not greater than (a + b), and $[C - (a + b)]^2$ cannot equal [2ab] when n = 3.

Therefore, there is no non-zero integer solution set for n = 3.

For n = 4 rewrite Equation 77 with n = 4.

(Note: The $[^n_n]$ term in Equation 77 is the third term when n = 4):

$$[C - (a + b)]^2 = a*b*([^4_2]C^{4-2}(2) - [^4_3]C^{4-3}(3a + 3b) + [^4_4]C^{4-4}(4a^2 + 6ab + 4b^2)) \div [C - (a + b)]^{4-2} \quad (82)$$

$$[C - (a + b)]^2 = a*b*(6C^2(2) - 4C(3a + 3b) + (4a^2 + 6ab + 4b^2)) \div [C - (a + b)]^2 \quad (83)$$

$$[C - (a + b)]^2 = 2ab[6C^2 - 6C(a + b) + (2a^2 + 3ab + 2b^2)] \div [C - (a + b)]^2 \quad (84)$$

Equation 84 is the same as Equation 59.

Let C = (a + b + x) and x = 1 in the right side of Equation 81 and solve for $[C - (a + b)]^2$.

$$[C - (a + b)]^2 = 2*a*b*(6C^2 - 6C(a + b) + (2a^2 + 3ab + 2b^2)) \quad (85)$$

$$[C - (a + b)]^2 = 2*a*b*(6(a + b + 1)^2 - 6(a + b + 1)(a + b) + (2a^2 + 3ab + 2b^2)) \quad (86)$$

$$[C - (a + b)]^2 = 2ab*(6(a + b + 1) + (2a^2 + 3ab + 2b^2)) \quad (87)$$

C must be greater than (a + b) for a common volume to exist (see Figures 1 and 2). The non-zero integer solution set for $[C - (a + b)]^2$ from Equation 9 is [2ab] when M is a non-zero integer and a and b are non-zero integer factors of $2M^2$ with $2M^2 = ab$ (see Equation 15). The quantity equal to $C - (a + b)]^2$ is always greater than [2ab] when C is greater than (a + b) when n = 4 from Equation 87.

Therefore, there is no non-zero integer solution set for n = 4.

For n = n repeat Equation 73:

$$[C - (a + b)]^n = ab([^n_2]C^{n-2}(2) - [^n_3]C^{n-3}(3a + 3b) + ...$$
$$\pm [^n_n] \{[^n_1]a^{n-2} + [^n_2]a^{n-3}b^1 + [^n_3]a^{n-4}b^2 + ... + [^n_{n-1}]b^{n-2}\}) \qquad (88)$$

The root of each side of Equation 88 (and 73) is [C - (a + b)].

Divide Equation 88 by $[C - (a + b)]^{n-2}$.

$$[C - (a + b)]^2 = a*b*([^n_2]C^{n-2}(2) - [^n_3]C^{n-3}(3a + 3b) + ...$$
$$\pm [^n_n] \{[^n_1]a^{n-2} + [^n_2]a^{n-3}b^1 + [^n_3]a^{n-4}b^2 + ... + [^n_{n-1}]b^{n-2}\}) \div [C - (a + b)]^{n-2} \qquad (89)$$

The non-zero integer solution set for $[C - (a + b)]^2$ is 2ab when the variables a and b are non-zero integer factors of $2M^2$ and M is a non-zero integer. Therefore, the quantity equal to $[C - (a + b)]^2$ in Equation 89 must equal 2ab for whole number solutions to exist. In addition, C must be greater than (a + b) for a common volume to exist (see Figures 1 and 2).

Note: The $[^n_n]$ in Equations 88 and 89 are binomial coefficients. The number of applicable terms and value in Equation 89 increase as n increases (i.e., number of applicable terms in Equation 89 equals n – 1). Completion of the nth power always begins with the $[^n_2]$ term in the series and ends with (n – 1) terms (see A in Figure 3). The quantity equal to $[C - (a + b)]^2$ is always greater than 2ab when n is greater than 2.

|  |  |  | Applicable terms - A |
|---|---|---|---|
| 1 |  |  |  |
| 1  1 |  |  |  |
| 1 -2  1 | n = 2 | A = 1 |  |
| 1 -3  3 -1 | n = 3 | A = 2 |  |
| 1 -4  6 -4  1 | n = 4 | A = 3 |  |
| 1 -5  10 -10  5 -1 | n = 5 | A = 4 |  |
| 1 -6  15 -20  15 -6  1 | n = 6 | A = 5 |  |
| 1 -7  21 -35  35 -21  7 -1 | n = 7 | A = 6 |  |
| 1 -8  28 -56  70 -56  28 -8  1 | n = 8 | A = 7 |  |
| 1 -9  36 -84  126 -126  84 -36  9 -1 | n = 9 | A = 8 |  |
| 1 -10  45 -120  210 -252  210 -120  45 -10  1 | n = 10 | A = 9 |  |

Figure 3. Binomial Coefficients for n = 2 to n = 10 for Equations 88 and 89. The highlighted numbers are the binomial coefficients for completing the $n^{th}$ power beginning with $[^n_2]$.

For example:

Completion of the square for n = 2 is $ab([^2_2]C^{2-2}(2)) = 2ab$

Completion of the cube for n = 3 is $ab([^3_2]C^{3-2}(2) - [^3_3]C^{3-3}(3a + 3b))$

   = $ab([3]C^1(2) - [1](3a + 3b)) = 3ab(2C - (a + b)$, the same as Equation 27.

Completion of the fourth power for n = 4 is $ab([^4_2]C^{4-2}(2) - [^4_3]C^{4-3}(3a + 3b) + [^4_4]C^{4-4}(2a^2 + 3ab + 2b^2))$

   = $2ab[6C^2 - 6C(a + b) + (2a^2 + 3ab + 2b^2)]$

The quantity $[C - (a + b)]^2$ in Equation 89 equals 2ab only when n = 2. Therefore, there is no non-zero integer solution set for any n greater than 2 since the quantity equal to $[C - (a + b)]^2$ is always greater than 2ab when n is greater than 2 and C is greater than (a + b).

Bibliography:

My cursory review of approaches to a proof of Fermat's Last Theorem did not find any similar to the one in this paper. Two references I relied on most are:

Fermat's Last Theorem (Vandiver 1914; Vandiver publisher; L'Intermédiaire des Math. 15, 217-218.; Swinnerton-Dwyer, P. Nature 364, 13--14, 1993) found at

http://mathworld.wolfram.com/FermatsLastTheorem.html.

Fermat's Enigma, by Simon Singh, 1997, Anchor Books, New York